\documentclass[a4paper,11pt]{article}
\usepackage{amsmath}
\usepackage{amsthm}
\usepackage[all]{xy}
\usepackage[english]{babel}
\usepackage[utf8]{inputenc}
\usepackage{graphicx}
\usepackage{xypic,color}
\usepackage{amsfonts, amssymb}

\newtheorem{thm}{Theorem}
\newtheorem{defn}[thm]{\textbf{Definition}}

\newtheorem{cor}[thm]{Corollary}

\newtheorem{prop}[thm]{Proposition}
\theoremstyle{definition}

\theoremstyle{definition}

\theoremstyle{remark}
\newtheorem{rem}[thm]{Remark}

\renewcommand{\k}{\Bbbk}

\newcommand{\de}{\Delta}
\newcommand{\al}{\alpha}

\newcommand{\ot}{\otimes}
\newcommand{\ma}{\mathcal}
\newcommand{\rt}{\rightarrow}

\title{Nearly Frobenius dimension on Frobenius algebras}
\author{Dalia Artenstein, Ana Gonz\'alez, Gustavo Mata.
}

\date{\today}

\begin{document}

	\maketitle

	\begin{abstract}
This article is divided into two parts.  In the first part we work over a field $\k$ and prove that the Frobenius space associated to a Frobenius algebra is generated as left A-module by the Frobenius coproduct. In particular, we prove that the Frobenius dimension coincides with the dimension of the algebra. In the second part we work with a commutative ring $k$. We introduce the concept of nearly Frobenius algebras in this context and construct solutions of the Yang-Baxter equation starting from elements in the Frobenius space. Also, we give a list of equivalent characterizations of nearly Frobenius algebras.
	
\end{abstract}

\section{Introduction}

\hspace{0.5cm}The concept of nearly Frobenius algebra is motivated by the result proved in \cite{CG}, which states that: the homology of the free loop space $H^{\ast}(LM)$ has the structure of a Frobenius algebra without counit. These objects were studied in \cite{GLSU} and their algebraic properties were developed in \cite{AGL15}, \cite{AGM1} and \cite{AGM2}. 
In \cite{AGL15} the notion of Frobenius space associated to a $\k$-algebra $A$, $\ma{E}_A$,  was defined and it was proved that if $A$ has finite dimension then $\ma{E}_A$ is a finite dimensional $\k$-vectorial space. The $\k$-dimension of $\ma{E}_A$ is denoted by $\operatorname{Frobdim}A$ and it follows that $\operatorname{Frobdim}A \leq (\operatorname{dim} A)^2$.
In \cite{AGM1} several classes of algebras with non trivial nearly Frobenius structure were described. Finally in \cite{AGM2} the notion of normalized nearly Frobenius structure was introduced and it was proved that an algebra $A$ has a normalized nearly Frobenius structure if and only if $A$ is a separable algebra if and only if the category of bimodules $ _A\mathcal{M}_A$ is semisimple.    \\

In this work we continue to deepen in the study of nearly Frobenius algebras. In the first part of the article we concentrate on answering the following question: What additional properties has the Frobenius space if we work with a Frobenius algebra? In the second part we go one step further in the generalization of the concept by extending the notion to algebras on commutative rings. Some properties are preserved and they are described in this work. Other characterizations of the concept emerge. Also a relationship with solutions of the quantum Yang-Baxter equation arises.\\

The article is organized as follows. In section $2$ all the algebras involved are of finite dimension over a field $\k$. We define an $A$-module structures on the space $\ma{E}_A$ and we obtain the following result:\\

{\bf{Theorem A}}
If $A$ is a Frobenius $\k$-algebra, then the Frobenius space $\mathcal{E}_A$ is generated, as left $A$-module, by  the Frobenius coproduct of $A$, let us call it $\de_0$. In other words, there exist $a\in A$ such that $ \de=a\star\de_0$ for all $\de\in\ma{E}_A$.\\

As a straight consequence of the previous theorem we have the next result for Frobenius algebras.\\

{\bf{Corollary A}}
If $A$ is a Frobenius algebra, then $$\operatorname{dim}_\k A=\operatorname{Frobdim}A.$$

In section $3$ the notion of nearly Frobenius algebra over a commutative rings is introduced. We construct solutions of the Yang-Baxter equation starting from elements in the Frobenius space. As a consequence of the results developed in this section we obtain the following characterizations of nearly Frobenius algebras:\\

{\bf{Corollary B}}
	Let $A$ be an algebra over a commutative ring $k$. Then the following statements are equivalent:
	\begin{enumerate}
		\item $A$ is a nearly Frobenius algebra.
		\item There exist $Q\in A\otimes_k A$ central.
		\item There exist $\eta_A:\operatorname{End}_k(A)\rt A\ot_kA$ an $A$-bimodule homorphism.
		\item The map $i: A\rt \ma{A}(R)$ given by $i(a)(b)=ab$ is a monomorphism of algebras, where $R\in\operatorname{End}_A(A\otimes A)$ is a homomorphism of right $A$-modules solution of the equation:
$$R^{13}\circ R^{12}=R^{23}\circ R^{13}=R^{12}\circ R^{23}.$$		
		
		\end{enumerate}

\section{Frobenius algebras over fields}
\subsection{Abrams Theorem}
In \cite{A99} Abrams proves a theorem that gives a characterization of  Frobenius algebras in terms of the existence of a coassociative counital comultiplication $\Delta:A\rightarrow A\otimes A$ which is a map of $A$-bimodules. In his construction, $\Delta(1)$ can be expressed from a basis and the dual one related to the bilinear form of the Frobenius algebras.  This expression will be useful to prove some new  results. Let us recall Abrams Theorem and some of its implications.
\begin{thm}\label{TM Abrams}
	An algebra $A$ is a Frobenius algebra if and only if it has a coassociative counital comultiplication $\Delta:A\rightarrow A\otimes A$ which is a map of $A$-bimodules.
\end{thm}
In the proof Abrams constructs the coproduct using the isomorphism of right (or left) $A$-modules $\lambda_r:A\to A^*$ ($\lambda_l:A\to A^*$) as follows 
\begin{equation}\label{eq1}
	\xymatrix{
A\ar[r]^(.4){\de}\ar[d]_{\lambda_r}& A\otimes A\\
A^*\ar[r]_(.4){\mu^*}&A^*\otimes A^*\ar[u]_{\lambda_r^{-1}\otimes\lambda_r^{-1}}
},	
\end{equation}
 where $\mu$ is the composition of the transposition $\tau:A\otimes A\rt A \ot A$ with the product $m$ of $A$.

Given a basis $\{e_1,\dots,e_n\}$ of the Frobenius algebra $A$ with bilinear form $B:A\times A\rt \k$, let $\Bigl\{e_1^\#,\dots, e_n^\#\Bigr\}$ be the \emph{dual basis} of $A$ relative to $B$, that means  the elements satisfy $B\Bigl(e_i^\#, e_j\Bigr)=\delta_{ij}$, for all $i,j=1,\dots,n$. On the other hand, $\{e_1^*,\dots, e_n^*\}$ denotes de basis of $A^*$ satisfying $e_i^*(e_j)=\delta_{ij}$. Then, using that $\lambda_r:A\rt A^* $  is given by $\lambda_r(a)(b)=B(a,b)$ and the fact that
$$\lambda_r\Bigl(e_i^\#\Bigr)(e_j)=B\Bigl(e_i^\#,e_j\Bigr)=\delta_{ij},\quad\mbox{for all}\quad j=1,\dots,n,$$
we conclude that $e_i^\#=\lambda_r^{-1}(e_i^*)$ for all $i=1,\dots, n$.

\begin{prop}\label{prop2}
	The coproduct $\Delta:A\rightarrow A\otimes A$ of the Frobenius algebra $A$ given in equation (\ref{eq1})  satisfies
	\begin{enumerate}
		\item[(i)] $\displaystyle{\de(1_A)=\sum_{i=1}^ne_i\otimes e_i^\#}$. 
		\item[(ii)] $\displaystyle{\de(x)=\sum_{i=1}^nxe_i\otimes e_i^\#=\sum_{i=1}^ne_i\otimes e_i^\#x}$, for all $x\in A$.
	\end{enumerate}
\end{prop}
\begin{proof}
	We write $\de(1)=\sum_ju_j\ot e_j^\#$. Then, using that $\de$ is a map of right $A$-modules and that $\varepsilon$, defined as $\varepsilon(x)=B(1_A,x)$, is the counit we have that  $$\sum_{j=1}^n\varepsilon\Bigl(e_j^\#x\Bigr)u_j=x\quad\mbox{for all $x\in A$}.$$
	In particular, if we take $\displaystyle{x=e_i \Rightarrow e_i=\sum_{j=1}^n\varepsilon\Bigl(e_j^\#e_i\Bigr)u_j=\sum_{j=1}^n\delta_{ij}u_j=u_i}.$ Then 
	 $$\de(1)=\sum_{j=1}^ne_j\ot e_j^\#.$$
	 The second claim is an immediate consequence of the fact that $\de$ is a map of $A$-bimodules.
\end{proof}

\subsection{Dimension of the Frobenius space on Frobenius algebras}
In this section we recall the definition of a nearly Frobenius algebra and the associated Frobenius space. Frobenius algebras are also nearly Frobenius algebras so we can compute their Frobenius space. We will prove that the dimension of the Frobenius space of a Frobenius algebra  coincides with the dimension of the algebra. \begin{defn}
	An algebra $A$ is a \emph{nearly Frobenius algebra} if it admits a map $\Delta:A\rt A\ot A$ of $A$-bimodules, that means, the following diagrams commute:
	$$\xymatrix{A\otimes A\ar[r]^{m}\ar[d]_{\Delta\otimes 1}& A\ar[d]^{\Delta}\\
		A\otimes A\otimes A\ar[r]_{1\otimes m}&A\otimes A
	},\quad \xymatrix{A\otimes A\ar[r]^{m}\ar[d]_{1\otimes\Delta}& A\ar[d]^{\Delta}\\
		A\otimes A\otimes A\ar[r]_{m\otimes 1}&A\otimes A
	}$$ 
\end{defn}

\begin{defn}
	The \emph{Frobenius space} associated to an algebra $A$ is the vector space of all the possible coproducts $\de$ that make it into a nearly Frobenius algebra ($\mathcal{E}_A$), see \cite{AGL15}. Its dimension over $\k$ is called the \emph{Frobenius dimension} of $A$, that is,
	$$\operatorname{Frobdim}A = \operatorname{dim}_\k\mathcal{E}_A.$$
\end{defn}

\begin{rem}

	The Frobenius space $\ma{E}_A$ associated to the algebra $A$ admits a structure of left $A$-module as follows
	$$\begin{array}{rcl}
		A\times \ma{E}_A&\rt& \ma{E}_A\\
		(a,\de)&\mapsto& a\ast\de
	\end{array}$$ $\bigl(a\ast\de\bigr)(1_A):=(1\otimes a)\de(1_A)=\sum x_i\otimes ay_i$, for all $a\in A$ where $\de(1_A)=\sum x_i\otimes y_i$.\\

Note that $a\ast\de\in\ma{E}_A$:
$$
\begin{array}{rcl}
(x\otimes 1)(a\ast\de)(1)&=&\displaystyle{(x\otimes 1)(1\otimes a)\de(1)=(1\otimes a)(x\otimes 1)\de(1)}\\
&=&(1\otimes a)\de(1)(1\otimes x)=(a\ast\de)(1)(1\otimes x).
\end{array}$$
It is clear that this operation defines a left action on $\ma{E}_A$.
\end{rem}
\begin{rem}
	If $A$ is a Frobenius algebra we can twist the previous action using the Nakayama automorphism associated to $A$ as follows: As $A$ is a Frobenius algebra there exists $\nu:A\rt A$ an automorphism of algebras such that $B(x,y)=B(y,\nu(x))$ for all $x,y\in A$. Then, we define the left action of $A$ on $\ma{E}_A$ as:
	$$(a\star\de)(1):=\Bigl(1\ot\nu^{-1}(a)\Bigr)\de(1)\quad\mbox{for all $a\in A$ and $\de\in\ma{E}_A$}.$$
	As $\nu$ is an algebra homomorphism so is $\nu^{-1}$. Then $\star$ is an action. \\
	
	Note that if $A$ is a symmetric algebra then the two actions previously defined coincide.
\end{rem}
\begin{prop}
	Let $A$ be a Frobenius algebra with Frobenius structure $(\de,\varepsilon)$ and $M=A\star\Delta\subset\mathcal{E}_A$ the left $A$-module generated by $\Delta$. Then, $$\operatorname{dim}_\k M=n.$$ 
\end{prop}
\begin{proof}
	We will prove that $\mathcal{D}=\bigl\{e_1\star\Delta,\dots,e_n\star\Delta\bigr\}$ is a basis of $M$ as vector space, where $\bigl\{e_1,\dots,e_n\bigr\}$ and $\Bigl\{e_1^\#,\dots, e_n^\#\Bigr\}$ are dual bases of $A$ related to the bilinear form $B$. \\
	First, note that $\mathcal{D}$ is a generator of $M$ by definition. Now, we will prove that $\mathcal{D}$ is LI: suppose that 
	$$\sum_{j=1}^n\al_j\bigl(e_j\star\Delta\bigr)=0$$ then, $$\sum_{j=1}^n\al_j\bigl(e_j\star\Delta\bigr)(1)=0.$$
	Remember that the Frobenius coproduct satisfies, by the Proposition \ref{prop2}, that $\displaystyle{\de(1)=\sum_{i=1}^ne_i\ot e_i^\#}$. Then
	$$\sum_{j=1}^n\al_j\sum_{i=1}^ne_i\ot \nu^{-1}\bigl(e_j\bigr)e_i^\#=0\Rightarrow \sum_{i,j=1}^n\alpha_je_i\ot \nu^{-1}\bigl(e_j\bigr)e_i^\#=0.$$
	Applying the map $(1\ot\varepsilon)$ we conclude that 
	$$0=\sum_{i,j=1}^n\alpha_je_i\varepsilon\Bigl(\nu^{-1}\bigl(e_j\Bigr)e_i^\#\bigr)=\sum_{i,j=1}^n\alpha_je_i\varepsilon\Bigl(e_i^\#e_j\Bigr)=\sum_{i=1}^n\alpha_ie_i\Rightarrow \alpha_i=0,\;\forall i=1,\dots,n.$$
	Therefore $\operatorname{dim}_\k M=n.$
\end{proof}

\begin{thm}\label{notsymmetric}
	If $A$ is a Frobenius $\k$-algebra, then the Frobenius space $\mathcal{E}_A$ is generated, as left $A$-module, by  the Frobenius coproduct of $A$, let us call it $\de_0$. In other words, there exist $a\in A$ such that $ \de=a\star\de_0$ for all $\de\in\ma{E}_A$.
\end{thm}
\begin{proof}
	We consider  $\bigl\{e_1,\dots,e_n\bigr\}$ and $\Bigl\{e_1^\#,\dots, e_n^\#\Bigr\}$ dual bases of $A$ related to the bilinear form $B$ and $\nu:A\rt A$ the Nakayama automorphism associated to $B$. Then, by Proposition \ref{prop2}, 
	$$\de_0(x)=\sum_{i=1}^nxe_i\ot e_i^\#$$
	is an $A$-bimodule homomorphism, and $\varepsilon(x)=B(1,x)=B(x,1)$ is the corresponding trace. \\
	Note that the dual bases satisfy that:
	$$x=\sum_{j=1}^n\varepsilon\Bigl(e_j^\#x\Bigr)e_j=\sum_{j=1}^n\varepsilon\bigl(xe_j\bigr)e_j^\#\quad\mbox{for all $x\in A$}.$$ 
	If we consider $\de\in\mathcal{E}_A$ we can write $\displaystyle{\de(1)=\sum_{i=1}^n e_i\ot b_i,}$ where $b_i\in A$ for $i=1,\dots,n$.\\
	Using that $\de$ is an $A$-bimodule homorphism and the previous identifications of all $x\in A$, we have that
	$$\sum_{i=1}^nxe_i\ot b_i=\sum_{i,j=1}^n\varepsilon\Bigl(e_j^\#xe_i\Bigr)e_j\ot b_i=\sum_{j=1}^ne_j\ot\left(\sum_{i=1}^n\varepsilon\Bigl(e_j^\#xe_i\Bigr)b_i\right)=\sum_{j=1}^ne_j\ot b_jx.$$
	Then $\displaystyle{b_jx=\sum_{i=1}^n\varepsilon\Bigl(e_j^\#xe_i\Bigr)b_i}$ for all $j=1,\dots,n$, $x\in A$. Now, we define a homomorphism $\sigma:A\rt A$ of $\k$-vector spaces by $\sigma\Bigl(e_j^\#\Bigr)=b_j$ for all $j=1,\dots, n$.\\
	Note that
	$$\sigma\Bigl(e_j^\#x\Bigr)=\sigma\left(\sum_{i=1}^n\varepsilon\Bigl(e_j^\#xe_i\Bigr)e_i^\#\right)=\sum_{i=1}^n\varepsilon\Bigl(e_j^\#xe_i\Bigr)\sigma\Bigl(e_i^\#\Bigr)=\sum_{i=1}^n\varepsilon\Bigl(e_j^\#xe_i\Bigr)b_i
	=b_jx.$$
	Then 
	$$\sigma\Bigl(e_j^\#x\Bigr)=b_jx,\quad\mbox{for all}\; j=1,\dots, n.$$ This proves that $\sigma $ is a homomorphism of right $A$-modules and so there exists an element $b\in A$ such that $\sigma(y)=b\cdot y$ for all $y\in A$. Finally we define $a=\nu(b)\in A$, then  $b_j=\sigma\Bigl(e_j^\#\Bigr)=\nu^{-1}(a)\cdot e_j^\#$, and $\de=a\star\de_0$.

\end{proof}	

\begin{cor}\label{symmetrc2}
	If $A$ is a Frobenius algebra, then $$\operatorname{dim}_\k A=\operatorname{Frobdim}A.$$
\end{cor}

\begin{rem}
	Note that the previous result in not true if the nearly Frobenius algebra is not Frobenius. An example of this remark is the path algebra $A=\k Q$, here
	$$Q=\xymatrix{
		\bullet_{1}\ar[r]^{\al}&\bullet_{2}}$$
	This algebra has, as a $\k$-vector space, dimension 3 and in \cite{AGL15} it is proven that $\operatorname{Frobdim}_\k A=1$. Then 
	$$\operatorname{Frobdim}_\k A<\operatorname{dim}_\k A.$$ 
	There are many ways to guarantee that this algebra is not  Frobenius. We will verify this by showing that no nearly coproduct can be completed to a Frobenuis one.\\
	The nearly Frobenius space is generated by the coproduct $\de$ defined as: 
	$$	\de(e_1)=\al\ot e_1,\quad \de(e_2)=e_2\ot\al,\quad \de(\al)=\al\otimes \al. $$
	Note that it is impossible to define a trace $\varepsilon:A\rt\k$ such that $(\varepsilon\ot 1)\de=\operatorname{Id}=(1\otimes\varepsilon)\de$.\\
	The next example shows that the Frobenius dimension can  also be greater than the dimension of the algebra.
	Let $A$ be the path algebra associated to the quiver
	$$Q=\xymatrix{
&\bullet_2\ar[rd]^\beta&\\
\bullet_{1}\ar[ur]^\alpha&&\bullet_3\ar[ll]^\gamma	
}$$ with the relation $\gamma\alpha=0$. \\
It is clear that this algebra has dimension $9$ as a $\k$-vector space. The Frobenius space has dimension $10$, and a general nearly Frobenius coproduct is expressed as 
$$\de(e_2)=a\beta\gamma\otimes\alpha,\quad \de(e_3)=a\gamma\otimes\alpha\beta$$
and   
$$\begin{array}{rcl}
\de(e_1)&=&a\bigl(e_1\otimes\alpha\beta\gamma+\alpha\beta\gamma\bigr)+b_1\alpha\otimes\gamma+b_2\alpha\otimes\beta\gamma+b_3\alpha\otimes\alpha\beta\gamma+b_4\alpha\beta\otimes\gamma\\
&=&b_5\alpha\beta\otimes\beta\gamma+b_6\alpha\beta\otimes\alpha\beta\gamma+b_7\alpha\beta\gamma\otimes\gamma+b_8\alpha\beta\gamma\otimes\beta\gamma+b_9\alpha\beta\gamma\otimes\alpha\beta\gamma
\end{array}$$
As before it is impossible to define a trace for these coproducts, which means that $A$ is not a Frobenius algebra.  
\end{rem}

\section{Frobenius structures over commutative rings}
In this section we will consider Frobenius and nearly Frobenius algebras over commutative rings instead of fields. In this context we will prove a version of Theorem \ref{notsymmetric}. In the last subsection we construct solutions of the Yang-Baxter equation starting from elements in the Frobenius space $\mathcal{E}_A$. Finally we give several characterizations to nearly Frobenius algebras over a commutative ring.

Throughout this section $k$ is a commutative ring.
\subsection{Frobenius algebras over commutative rings}

\begin{defn}
	A $k$-algebra $A$ is \emph{Frobenius} if $A$ is finitely generated and projective as a $k$-module, and there exists a left $A$-module isomorphism 
	$$\lambda_l:A\rightarrow A^*.$$
\end{defn}

\begin{rem}\label{rem12}
	Take an arbitrary ring extension $B\subset A$ and let
	$$F: {}_A\mathcal{M}\rt {}_B\mathcal{M}$$
	be the restriction of scalars functor. Then $F$ has a left adjoint, the induction functor
	$$G=A\otimes_B-: {}_B\mathcal{M}\rt {}_A\mathcal{M}$$	
	and a right adjoint, the coinduction functor
	$$T=Hom_B(A,-):  {}_B\mathcal{M}\rt {}_A\mathcal{M}.$$
	The left $A$-module structure on $T(N)=Hom_B(A,N)$ is given by $(a\cdot\gamma)(b)=\gamma(ba)$, for all $a,b\in A$ and $\gamma\in Hom_B(A,N)$. For a left $B$-module map $f$, we put $T(f)=f\circ-$.
\end{rem}

\begin{defn}
	A covariant functor $F:\ma{C}\rt\ma{D}$ is called a \emph{ Frobenius functor} if it has isomorphic left and right adjoints.
\end{defn}

\begin{thm}[Theorem 2.4 of \cite{CIM}]\label{t14}
	Let $B\subset A$ be an arbitrary ring extension. Then the following statements are equivalent:
	\begin{enumerate}
		\item The restriction of scalars functor $F: {}_A\mathcal{M}\rt {}_B\mathcal{M}$ is Frobenius (or, equivalently, $B\subset A$ is a Frobenius extension).
		\item There exist a pair $(\de,\varepsilon)$ such that
		\begin{enumerate}
			\item[(a)] $\de:A\rt A\ot_BA$ is an $A$-bimodule map,
			\item[(b)] $\varepsilon:A\rt B$ is a $B$-bimodule map, and a counit for $\de$, that is $$(1\otimes_B\varepsilon)\de=(\varepsilon\ot_B1)\de=\operatorname{Id}_A.$$
		\end{enumerate}
	\end{enumerate}
	\end{thm}

\begin{cor}[Corollary 2.6 of \cite{CIM}]
	For an algebra $A$ over a commutative ring $k$, the following statements are equivalent:
\begin{enumerate}
	\item $A$ is a Frobenius algebra.
	\item There exist a coalgebra structure $\bigl(A,\de,\varepsilon\bigr)$ on $A$ such that the coproduct $\de:A\rt A\ot_kA$ is an $A$-bimodule map.
	\item  There exist $Q=\sum_{i=1}^n e_i\ot e_i^\#\in A\ot_kA $ and $\varepsilon\in A^*$ such that $Q$ is $A$-central and the normalizing Frobenius condition
	$$\sum_{i=1}^n\varepsilon\bigl(e_i\bigr)e_i^\#=\sum_{i=1}^n e_i\varepsilon\Bigl(e_i^\#\Bigr)=1_A$$
	is satisfied. $(Q,\varepsilon)$ is called a Frobenius pair.
\end{enumerate}
	\end{cor}

\begin{rem}
	As a consequence of the previous results we have the following properties:
	\begin{enumerate}
		\item[(1)] $Q$ is $A$-central then $\displaystyle{\sum_{i=1}^n xe_i\ot e_i^\#=\sum_{i=1}^n e_i\ot e_i^\#x\quad\mbox{for all $x\in A$}}$.
		\item[(2)] The normalizing condition implies that: 
		$\displaystyle{x=\sum_{i=1}^n\varepsilon\bigl(xe_i\bigr)e_i^\#=\sum_{i=1}^n\varepsilon\Bigl(e_i^\#x\Bigr)e_i}$\; for all $x\in A.$
		\item[(3)] The sets $\bigl\{e_1,\dots, e_n\bigr\}$ and $\Bigl\{e_1^\#,\dots, e_n^\#\Bigr\}$ are generators of $A$ as $k$-module. We call them \emph{dual generators} of $A$. 
		\item[(4)] The element $Q$ does not depend on the choice of dual generators $e_i$, $e_i^\#$, $i=1,\dots,n$ of $A$.
	\end{enumerate}
	
\end{rem}

\subsection{Nearly Frobenius algebras over commutative rings}
\begin{defn}
	An associative $k$-algebra $A$ is a \emph{nearly Frobenius algebra} if it admits an homomorphism $\Delta:A\rt A\ot_k A$ of $k$-modules such that $\de$ is an homomorphism of $A$-bimodules.

	Let $(A,\Delta_A)$ and $(B,\Delta_B)$ be two nearly Frobenius algebras. An homomorphism $f:A\rightarrow B$ is a \emph{nearly Frobenius homomorphism} if it is a morphism of $k$-algebras and the following diagram commutes.
	$$\xymatrix{A\ar[r]^{f}\ar[d]_{\Delta_A}& B\ar[d]^{\Delta_B}\\
		A\otimes_k A\ar[r]_{f\otimes f}&B\otimes_k B
	}.$$
\end{defn}

\begin{rem}
	If we have an associative $k$-algebra $A$ we can consider the set $\mathcal{E}_A$ of all coproducts $(\de:A\rt A\otimes_k A)$ that make it into a nearly Frobenius algebra. The set $\mathcal{E}_A$ is a $k$-module. The proof that $(\mathcal{E}_A,+)$ is an abelian group is found in \cite{AGL15}. We define the operation $\cdot: k\times\mathcal{E}_A\rt\mathcal{E}_A$ as follows: $$(r\cdot\de)(x)=(r\otimes_k 1)\de(x)=\sum r\cdot x_1\otimes x_2,\; \mbox{if}\; \de(x)=\sum \cdot x_1\otimes_k x_2,$$ for all $r\in k$ and $\de\in\mathcal{E}_A$. It is easy to see that this operation defines an action of $A$ on $\ma{E}_A$.
	
	Finally, we define the \emph{Frobenius space over $k$} associated to an algebra $A$ as the $k$-module $\mathcal{E}_A$. \\ Note that $\mathcal{E}_A=Hom_{A\mbox{-}A}(A,A\otimes_k A)$.
\end{rem}

\begin{defn} Let $A$ be an associative $k$-algebra. 
	Denote by $C_{A\otimes_k A}(A)$ the set of $A$-centralising elements of $A\otimes_k A$, i.e
	$$C_{A\otimes_k A}(A)=\left\{\sum_ix_i\otimes_k y_i\in A\otimes_k A:\;\sum_iax_i\otimes_k y_i=\sum_ix_i\otimes_k y_ia,\;\mbox{for all $ a\in A$}\right\}.$$
\end{defn}

\begin{rem}
	The evaluation homomorphism
	$$\Psi:\mathcal{E}_A\rt C_{A\otimes_k A}(A)\quad\mbox{with}\quad\de\mapsto\de(1_A)$$
	is an isomorphism of $k$-modules.\\
	In the previous remark we see that $\mathcal{E}_A$ is an $k$-module. The set $C_{A\otimes_kA}(A)$ is a natural $k$-submodule of $A\otimes_kA$ and it is clear that $\Psi$ is a morphism of $k$-modules. \\
	Note that if we know the value of $\de$ on $1_A$ then we can determine the value on all element of $A$ as: $\de(x)=x\cdot\de(1_A)$. Then $\Psi$ is surjective. In a similar way, if we have two coproducts $\de_1$ and $\de_2$ such that $\de_1(1_A)=\de_2(1_A)$, then $\de_1=\de_2$ because they are morphisms of $A$-bimodules. Then $\Psi$ is injective.
\end{rem}

\begin{rem}
	If $A$ is a Frobenius algebra the space $C_{A\otimes_k A}(A)$ has structure of left $A$-module as follows 
	$$a\star \sum_ix_i\otimes_k y_i:=\sum_ix_i\otimes_k \nu^{-1}(a)y_i\quad \mbox{for all $a\in A$},$$
	where $\nu:A\rt A$ is the Nakayama automorphism of $A$.
\end{rem}
\begin{thm}
	Let $A$ be a Frobenius algebra over a commutative ring $k$ with Nakayama automorphism $\nu:A\rt A$ and let $\bigl\{e_1,\dots, e_n\bigr\}$, $\Bigl\{e_1^\#,\dots, e_n^\#\Bigr\}$ be dual generators of $A$ and 	$\displaystyle{Q=\sum_{i=1}^ne_i\otimes e_i^\#\in A\ot_k A}$. Then 
	$C_{A\otimes_kA}(A)$ is a free rank one left $A$-submodule of $A\otimes_kA$ with basis $\bigl\{Q\bigr\}$.
\end{thm}
\begin{proof}
	It is clear that $C_{A\otimes_kA}(A)$ is a left $A$-submodule of $A\otimes_kA$ with the action $\star.$\\
	Now consider the left $A$-module $A\star Q$ and note that $A\star Q\subset C_{A\otimes_kA}(A)$. First we prove that $A\star Q$ is free: suppose that $a\star Q=0$ then
	$$0=(1\otimes\varepsilon)(a\star Q)=(1\otimes\varepsilon)\left(\sum_ie_i\otimes \nu^{-1}(a)e_i^\#\right)=\sum_i\varepsilon\bigl(\nu^{-1}(a)e_i^\#\bigr)e_i=\sum_i\varepsilon\bigl(e_i^\#a\bigr)e_i=a,$$
	then $a=0$. To finish it is enough to prove that $A\star Q=C_{A\otimes_kA}(A)$.\\
	Using localization to maximal ideals of $k$ arguments we can suppose that $A$ is a free $k$-module and reproduce the prove of Theorem \ref{notsymmetric} to guarantee  that $A\star Q=C_{A\otimes_kA}(A)$.
	
\end{proof}

The next result is an adaptation of Theorem \ref{t14} to the case of nearly Frobenius algebras.
\begin{thm}\label{extension}
	Let $B\subset A$ be an arbitrary ring extension. Then the following statements are equivalent:
	\begin{enumerate}
		\item[(a)] There exists an $A$-bimodule map $\de:A\rt A\otimes_BA$ .
		\item[(b)] There exists a natural transformation $\eta:T\rt G$, with $G, T: {}_B\mathcal{M}\rt {}_A\mathcal{M}$ defined in Remark \ref{rem12}.
	\end{enumerate}
\end{thm}
\begin{proof}
	(a)$\Rightarrow$(b) For $N\in {}_B\mathcal{M}$ we define $\eta_N:Hom_B(A,N)\rt A\otimes_BN$ as\\ $\eta_N(f)=\sum e_1\otimes_B f(e_2)$, where $\de(1_A)=\sum e_1\otimes_B e_2$.\\
	Let us verify that $\eta_N$ is an $A$-module morphism: $$a\cdot\eta_N(f)=a\cdot\left(\sum e_1\otimes_B f\bigl(e_2\bigr)\right)=\sum ae_1\otimes_B f\bigl(e_2\bigr).$$
	On the other hand, $$\eta_N(a\cdot f)=\sum e_1\otimes_B(a\cdot f)\bigl(e_2\bigr)=\sum e_1\otimes_Bf\bigl(e_2a\bigr).$$
	Using that $\de\in\mathcal{E}_A$ we have that $$(1\otimes f)\left(\sum ae_1\otimes_B e_2\right)=(1\otimes f)\left(\sum e_1\otimes_B e_2a\right)$$
	Then, $\sum ae_1\otimes_B f\bigl(e_2\bigr)=\sum e_1\otimes_Bf\bigl(e_2a\bigr)$, therefore $a\cdot\eta_N(f)=\eta_N(a\cdot f)$.\\
	Let us prove that $\eta$ natural transformation: For $M, N\in{}_B\mathcal{M}$ and $f\in Hom_B(M,N)$ 
	the following diagram must commute
	$$\xymatrix{Hom_B(A,M)\ar[r]^{\eta_M}\ar[d]_{Hom_B(A,f)}&A\otimes_BM\ar[d]^{A\otimes_Bf}\\
		Hom_B(A,N)\ar[r]_{\eta_N}&A\otimes_BN}$$ 
	Let be $\varphi\in Hom_B(A,M)$, then $(A\otimes_Bf)\eta_M(\varphi)=\sum e_1\otimes_Bf\varphi(e_2)=\eta_N(f\circ\varphi)=\eta_N\circ Hom_B(A,f)(\varphi)$.\\
	(b)$\Rightarrow$(a) For $a\in A$, we define $\varphi_a: A\rt A$ by $\varphi_a(b) = ba$. Then the map $\varphi: A\rt Hom_B(A,A)$ mapping $a$ to $\varphi_a$ is a left $A$-module morphism. Now take $$\de=\eta_A(\varphi)\quad (i.e. \de(a)=\eta_A(\varphi_a)).$$
	The map $\de$ is left $A$-module morphism, because $\eta_A$ and $\varphi$ are left A-module morphisms.\\
	As $\eta$ is a natural transformation we have the following commutative diagram
	$$\xymatrix{
		Hom_B(A,A)\ar[r]^{\eta_A}\ar[d]_{Hom_B(A,\varphi_a)}&A\otimes_BA\ar[d]^{A\otimes_B\varphi_a}\\
		Hom_B(A,A)\ar[r]_{\eta_A}&A\otimes_BA
	}$$
	Applying the diagram to $I_A =\varphi_1$, we have that
	$$\sum e_1\otimes_B e_2a=\bigl(A\otimes_B\varphi_a\bigr)\de(1_A)=\bigl(A\otimes_B\varphi_a\bigr)\circ\eta_A\bigl(\varphi_1\bigr)=\eta_A\circ Hom_B(A,\varphi_a)(\varphi_1)$$$$=\eta_A\bigl(\varphi_a\circ\varphi_1\bigr)=\eta_A(\varphi_a)=\de(a).$$
	Then $\de$ is a right $A$-module morphism.
\end{proof}

\subsection{The quantum Yang-Baxter equation}
The quantum Yang-Baxter equation (QYBE) is one of the basic equations in mathematical physics that lies in the foundation of the theory of quantum groups. In the last years, many solutions of this equation were found, and the related algebraic structures (Hopf algebras) have been intensively studied.

In \cite{BFS} it is shown that every Frobenius algebra over a commutative ring determines a class of solutions of the quantum Yang-Baxter equation. In this article we adapt the previous result to nearly Frobenius algebras. 

As before $k$ is an associative commutative ring with $1$ and $A$ an associative $k$-algebra. Let $\mathrm{R}\in\operatorname{End}(A\otimes_kA)$, $\mathrm{R}^{12}=\mathrm{R}\otimes\operatorname{Id}$, $\mathrm{R}^{23}=\operatorname{Id}\otimes \mathrm{R},$ and $\mathrm{R}^{13}=(\operatorname{Id}\otimes\tau)(\mathrm{R}\otimes\operatorname{Id})(\operatorname{Id}\otimes \tau)$ maps in $\operatorname{End}(A\otimes_kA\otimes_kA)$, where $\tau:A\otimes_kA\rt A\otimes_kA$ is the twist mapping. We shall say that $\mathrm{R}$ is a \emph{solution of the quantum Yang-Baxter equation} (QYBE), if
$$\mathrm{R}^{12}\circ \mathrm{R}^{13}\circ \mathrm{R}^{23}= \mathrm{R}^{23}\circ\mathrm{R}^{13}\circ\mathrm{R}^{12}.$$

\begin{prop}
Let $A$ be a $k$-algebra and $Q\in C_{A\otimes_kA}(A)\simeq\ma{E}_A$ (the Frobenius space of $A$). Then 
\begin{enumerate}
	\item $Q^{13}Q^{12}=Q^{23}Q^{13}=Q^{12}Q^{23}$,
	\item $Q^{12}Q^{23}Q^{12}=Q^{23}Q^{12}Q^{23}$,
\end{enumerate}
where $Q^{12}=\sum_ix_i\otimes y_i\otimes 1$, $Q^{13}=\sum_ix_i\otimes 1\otimes y_i$ and $Q^{23}=\sum_i1\otimes x_i\otimes y_i$ if $Q=\sum_ix_i\otimes y_i$.
\end{prop}
\begin{proof}
\begin{enumerate}
	\item $$Q^{13}Q^{12}=\left(\sum_ix_i\otimes 1\otimes y_i\right)\left(\sum_jx_j\otimes y_j\otimes 1\right)=\sum_i(x_i\otimes 1)\sum_{j}(x_j\otimes y_j)\otimes y_i$$
	$$=\sum_{i,j}x_j\otimes y_jx_i\otimes y_i=\sum_jx_j\otimes(y_j\otimes 1)\sum_ix_i\otimes y_i=\sum_{i,j}x_j\otimes x_i\otimes y_iy_j=Q^{23}Q^{13},$$
	$$Q^{12}Q^{23}=\left(\sum_ix_i\otimes y_i\otimes 1\right)\left(\sum_j1\otimes x_j\otimes y_j\right)=\sum_ix_i\otimes (y_i\otimes 1)\sum_jx_j\otimes y_j$$
	$$=\sum_ix_i\otimes\sum_jx_j\otimes y_jy_i=\left(\sum_j1\otimes x_j\otimes y_j\right)\left( \sum_ix_i\otimes 1\otimes y_i\right)=Q^{23}Q^{13}.$$
	\item $$\bigl(Q^{12}Q^{23}\bigr)Q^{12}=\bigl(Q^{23}Q^{13}\bigr)Q^{12}=Q^{23}\bigl(Q^{13}Q^{12}\bigr)=Q^{23}\bigl(Q^{12}Q^{23}\bigr).$$
\end{enumerate}	
\end{proof}

Using the previous proposition and Lemma 2. of \cite{BFS} we have the following result.  
\begin{cor}
	Let $\mathrm{R}=Q\tau$, where $\tau:A\otimes_kA\rt A\otimes_kA$ is the twist map. Then $\mathrm{R}$ is a solution of the QYBE.
\end{cor}

The following results are based on the article \cite{CIM}.
\begin{defn}
	Let $A$ be an algebra over a commutative ring $k$ and $R\in\operatorname{End}_k(A\otimes_kA)$ a solution of the equation
	\begin{equation}\label{equation2}
		R^{13}\circ R^{12}=R^{23}\circ R^{13}=R^{12}\circ R^{23}.
	\end{equation}
We define the $k$-subalgebra $\ma{A}(R)$ of $\operatorname{End}_k(A)$ as 
$$\ma{A}(R):=\bigl\{f\in \operatorname{End}_k(A): f\cdot R=R\cdot f\bigr\},$$
with unit $\operatorname{Id}_A$.
\end{defn}
\begin{rem}
	Note that $\ma{A}(R)=\bigl\{f\in \operatorname{End}_k(A): (f\otimes\operatorname{Id}_A)\circ R=R\circ \bigl(\operatorname{Id}_A\otimes f\bigr)\bigr\}$ and $R$ is $\ma{A}(R)$-central. Moreover, if $\operatorname{End}_k(A)$ is flat as a $k$-module
	$R\in\ma{A}(R)\otimes\ma{A}(R)$: \\
	$\ma{A}(R)$ can be seen as  $\operatorname{Ker}(\varphi)$, with $\varphi:\operatorname{End}_k(A)\rt\operatorname{End}_k(A\otimes A) $ defined by
	$$\varphi(f)=(f\ot\operatorname{Id}_A)\circ R-R\circ(\operatorname{Id}_A\otimes f).$$
	Assuming $\operatorname{End}_k(A)$ is flat as a $k$-module, $\ma{A}(R)\otimes \operatorname{End}_k(A)=\operatorname{Ker}(\varphi\ot \operatorname{Id}_{\operatorname{End}_kA})$. In particular, if $R=\sum R^1\ot R^2$
	$$\begin{array}{rcl}
	\bigl(\varphi\ot \operatorname{Id}_{\operatorname{End}_kA}\bigr)(R)&=&\sum R^1\circ R^1\ot R^2\ot R^2 -\sum R^1\ot R^2 \circ R^1\ot R^2\\
	&=&R^{13}\circ R^{12}-R^{12}\circ R^{23}=0
	\end{array}$$
Then $R\in\ma{A}(R)\otimes \operatorname{End}_k(A)$. In similar way, using the other condition we can see that $R\in\operatorname{End}_k(A)\ot\ma{A}(R).$ Then  $R\in\ma{A}(R)\otimes\ma{A}(R)$.
\end{rem}

\begin{thm}\label{monomorphism}
	For an algebra $A$ over a commutative ring $k$, the following statements are equivalent:\begin{enumerate}
		\item $A$ is a nearly Frobenius algebra.
		\item The map $i: A\rt \ma{A}(R)$ given by $i(a)(b)=ab$ is a monomorphism of algebras, where $R\in\operatorname{End}_A(A\otimes A)$ is a homomorphism of right $A$-modules solution of (\ref{equation2}). 
	\end{enumerate}
\end{thm}
\begin{proof}
	$(1)\Rightarrow (2)$ Let $\de$ be the nearly Frobenius coproduct and $Q=\de(1)=\sum_ie_i\otimes e^i$ the corresponding $A$-central element. Then, the map 
	$$R:A\otimes A\rt A\otimes A,\quad R(a\otimes b)=\sum_ie_ia\ot e^ib$$
	is a solution of (\ref{equation2}). Consequently, we can construct the algebra $\ma{A}(R)=\bigl\{f\in\operatorname{End}_k(A):\bigl(f\otimes \operatorname{Id}_A\bigr)\circ R=R\circ\bigl(\operatorname{Id}_A\ot f\bigr)\bigr\}\subseteq\operatorname{End}_k(A)$,  and define the injection map $i:A\rt\operatorname{End}_k(A)$ as $i(a)(b)=ab$, for all $a,b\in A$. We will prove that $i(A)\subset\ma{A}(R)$.
	$$\begin{array}{rcl}
			\bigl(i(a)\ot \operatorname{Id}_A\bigr)R(x\otimes y)&=&\displaystyle{\bigl(i(a)\ot \operatorname{Id}_A\bigr)\bigl(\sum_ie_ix\otimes e^iy\bigr)=\sum_iae_ix\ot e^iy}\\
			&=&\displaystyle{\bigl(\sum_iae_i\ot e^i\bigr)(x\ot y)=\bigl(\sum_ie_i\ot e^ia\bigr)(x\ot y)}\\
			&=&\displaystyle{\sum_ie_ix\ot e^iay=\bigl(\sum_ie_i\ot e^i\bigr)\bigl(x\ot ay\bigr)}\\
			&=&\displaystyle{R(x\otimes ay)=R\bigl(\operatorname{Id}_A\ot i(a)\bigr)(x\ot y),\;\forall x,y\in A.}
	\end{array}
$$ Then $	\bigl(i(a)\ot \operatorname{Id}_A\bigr)\circ R=R\circ\bigl(\operatorname{Id}_A\ot i(a)\bigr)$, and $i(A)\subset\ma{A}(R)$.\\
Note that $R$ is homomorphism of right $A$-module by definition: $R\bigl((a\otimes b)\cdot x\bigr)=R(a\otimes bx)=Q(a\otimes bx)=Q(a\otimes b)(1\otimes x)=R(a\otimes b)(1\otimes x)=R(a\otimes b)\cdot x$, for all $a,b,x\in A$.

$(2)\Rightarrow (1)$ We define $Q:=R(1\otimes 1)=\sum_ie_i\otimes f_i$. By Corollary  24 it is enough to prove that $Q$ is $A$-central.\\
As $i(a)\in\ma{A}(R)$ we have that $i(a)\cdot R=R\cdot i(a)$. Then $$\bigl[(i(a)\otimes \operatorname{Id})\circ R\bigr](1\otimes 1)=\bigl[R\cdot (\operatorname{Id}\otimes i(a))\bigr](1\otimes 1)\Rightarrow$$
$$\sum_ii(a)(e_i)\otimes f_i=\sum_iae_i\otimes f_i=R(1\otimes a)=R\bigl((1\otimes 1)\cdot a)\bigr)=R(1\otimes 1)\cdot a=Q\cdot a.$$
Then $a\cdot Q=Q\cdot a,$ for all $a\in A$.
\end{proof}

Finally we summarize some equivalences to nearly Frobenius algebras over a commutative ring. 

\begin{cor}
	Let $A$ be an algebra over a commutative ring $k$. Then the following statements are equivalent:
	\begin{enumerate}
		\item $A$ is a nearly Frobenius algebra.
		\item There exist $Q\in A\otimes_k A$ central.
		\item There exist $\eta_A:\operatorname{End}_k(A)\rt A\ot_kA$ an $A$-bimodule homorphism.
		\item The map $i: A\hookrightarrow \ma{A}(R)$ given by $i(a)(b)=ab$ is a monomorphism of algebras, where $R\in\operatorname{End}_A(A\otimes A)$ is a homomorphism of right $A$-modules solution of (\ref{equation2}).
		\end{enumerate}
\end{cor}
\begin{proof}
$(1)\Leftrightarrow (2) \Leftrightarrow (3)$ follows from Theorem \ref{extension} and $(1) \Leftrightarrow (4)$ follows from Theorem \ref{monomorphism}.\end{proof}

\bibliographystyle{elsarticle-harv}

\end{document}